\documentclass[12pt]{article}
\usepackage{amsmath}
\usepackage{amsfonts}
\usepackage{amssymb}

\usepackage{graphicx}%
\usepackage{multirow}%
\usepackage{mathrsfs}%
\usepackage[title]{appendix}%
\usepackage{xcolor}%
\usepackage{textcomp}%
\usepackage{manyfoot}%
\usepackage{booktabs}%
\usepackage{algorithm}%
\usepackage{algorithmicx}%
\usepackage{algpseudocode}%
\usepackage{listings}%
\title{Limit groups with respect to Thompson's group $F$ and other hereditarily separating groups}
\author{A. Ivanov and R. Zarzycki}

\newtheorem{theorem}{Theorem}[section]
\newtheorem{proposition}[theorem]{Proposition}
\newtheorem{corollary}[theorem]{Corollary}

\newtheorem{fact}[theorem]{Fact}

\newtheorem{theo}{Theorem}

\newtheorem{definition}[theo]{Definition}

\newtheorem{theore}{Theorem} 
\newtheorem{remark}[theore]{Remark}
\newtheorem{example}[theore]{Example}

\begin{document} 
\topmargin = 12pt
\textheight = 630pt 
\footskip = 39pt 

\maketitle

\begin{quote}
{\bf Abstract} 
We prove that a limit group over Thompson's group $F$ cannot be an HNN extension of $F$ with finitely generated associated subgroups. 
On the other hand, we give an example of an $F$-limit  group which is a centralized HNN extentions of $F$. 
We also study relatively limit groups with respect to $F$ and some other known groups of homeomorphisms of topological spaces.

\bigskip 


\end{quote}




\section{Introduction}\label{sec1}

 The notion of \emph{limit groups} was introduced by Sela in his work on the characterization
of the elementary equivalence of free groups \cite{Sel}. 
The idea has been extended in the paper by Champetier and Guirardel \cite{GC}, 
where the authors look at limit groups as limits of convergent sequences in spaces of \emph{marked groups}. 
Since then, this approach has grown into a significant branch of group theory. 
Among limits over the class of free groups, limits over other prominent groups are studied; 
for example, see papers of 
Guyot and Stalder  concerning the classes of dihedral and Baumslag-Solitar groups, \cite{G}, \cite{SG}, the paper of  Casals-Ruiz, Duncan and Kazachkov concerning right-angled Artin groups, \cite{CDK}, or the paper of Kochloukova and Zalesskii on pro-p analogs of limit groups, \cite{KZ}. 

Since Thompson's group $F$ has remained one of the most interesting objects in geometric group
theory, it is natural to consider limits with respect to $F$. 
This work was initiated by Akhmedov, Stain and Taback in \cite{Ast}, Brin in \cite{Bri} and Zarzycki in \cite{Zar}. 
We especially mention Brin's theorem that the free group $\mathbb{F}_2$ is an $F$-limit group, \cite{Bri}, 
which shows that Sela's limit groups are also $F$-limit ones. 

Our main results concern the question of when free constructions over $F$ (or some other related groups) are limit groups over $F$. 
In particular, in Section 2 we show that some HNN extensions of $F$ can be realized as limits of $F$. 
On the other hand, in Section 3 we prove: 
\begin{quote} 
\begin{itemize} 
\item $F$ cannot be a free summand in an $F$-limit group; 
\item any HNN extension of $F$ with finitely generated associated subgroups is not a limit over $F$. 
\end{itemize} 
\end{quote} 
These results strongly contrast with the situation of limits of free groups. 
In the latter case, a limit group is obtained from free groups by an iteration of separated amalgamated free
products or separated HNN extensions. 

 In our study, we employ mainly two technical tools. 
The first is taken from the recent paper of the authors \cite{Iva-Za} where solutions of systems of inequalities are studied 
in  groups having \emph{hereditarily separating actions}.
The second tool is the result of the second author on the existence of mixed identities in Thompson's group $F$; see \cite{SZ}. 
(To a large extent, the present paper can be viewed as the next step of the project started in 
\cite{Zar} and continued in \cite{SZ} and \cite{Iva-Za}.) 

Our methods can be used for some other well-known groups. 
It turns out that groups having hereditarily separating actions arise frequently in some areas of topological groups.
Furthermore, mixed identities of such groups have become an intensively studied topic in group 
theory; 
see \cite{HO}, \cite{Os} or \cite{BodSnTh}. 
Slightly extending the context, these examples can be viewed as natural objects of the theory of limit groups. 
We contribute in these respects by applying our approach to branch groups and some highly homogeneous groups.   

The structure of the paper is as follows. 
In the rest of this introduction we review some material on limit groups, mixed identities and Thompson's group $F$. 
Additional motivating comments are given there. 
In Section 2 we study examples of limit groups over $F$. 
Here we strongly use \cite{Iva-Za}. 
Necessary preliminaries will be given in this section. 
In Section 3 we use mixed identities to prove that any $F$-limit group 
is not an HNN extension with finitely generated associated subgroups. 
Section 4 is devoted to other groups with hereditarily separating actions.

\subsection{Marked groups and limit groups}

 A \emph{marked group} $(G,S)$ is a group G with a distinguished set of generators 
$S = (s_{1}, s_{2}, \ldots , s_{m})$. 
For a fixed $m$, let $\mathcal{G} _{m}$ be the set of all
$m$-generated groups marked by $m$ generators (up to isomorphism of marked groups).
Following \cite{GC} we put a certain metric on $\mathcal{G}_{m}$:  two marked
groups $(G,S), (G',S')\in\mathcal{G} _{m}$ are at distance less or equal to $e^{-r}$ if they
have exactly the same relations of length at most $r$. 
The set $\mathcal{G} _{m}$ equipped
with this metric is a compact space \cite{GC}. \emph{Limit groups} are simply limits of
convergent sequences in this metric space.
The following notion is basic in the paper. 

\begin{definition} \label{dlg} 
Let $G$ be an $m$-generated group. 
A marked group in $\mathcal{G} _{m}$ is
a limit group over $G$ (shortly $G$-limit) if it is a limit of marked groups each
isomorphic to $G$.
\end{definition}
Quite often in the literature, a slightly weaker definition is considered, 
where marked groups that appear are assumed to be isomorphic to subgroups of $G$. 
In our paper, when limit groups with respect to Thompson's group $F$ are considered, 
we usually take a sequence, $\{ (g_{1,n}, \ldots , g_{m,n} ) \, | \, n\in \mathbb{N}\}$ of elements from 
$F$ and the corresponding sequence of groups marked by $m+2$ elements, 
$G_{n}=(F,(\mathsf{x}_0,\mathsf{x}_1,g_{1,n},\ldots ,g_{m,n}))$, $n\in\mathbb{N}$, where $\mathsf{x}_0$ and $\mathsf{x}_1$
are the standard generators of $F$.

 It has been shown in \cite{GC}, \cite{Sel}  and \cite{FGM} that in the case of free groups some standard constructions
can be obtained as limits of free groups. 
For example, the following groups are the limits of convergent sequences in the space of free groups marked by three
elements: the free group of rank $3$, the free abelian group of rank $3$ or an HNN-extension
over a cyclic subgroup of the free group of rank $2$. 
All non-exceptional surface groups (see \cite{BaB}, \cite{BaG}) form another broad class of interesting examples.  
It is proved in \cite{Sel} and \cite{GC} that a finitely generated group is a limit of free groups if and only if it is fully residually free.

\paragraph{$\bullet$ Limit groups with infinitely many markers} 

In order to apply our method to infinitely generated groups (for example, oligomorphic ones; 
see \cite{BodSnTh}), we consider the following situation. 
Let $G$ be a group, $\bar{x}$ be a finite tuple of variables, and $W(\bar{x})$ be a subset of $G*\mathbb{F}_{\bar{x}}$.    
A set $p(\bar{x})$ of relations of the form $w(\bar{x}) = 1$ or $w(\bar{x}) \not=1$ where $w(\bar{x}) \in W(\bar{x})$, 
is called an {\em atomic type over} $W(\bar{x})$ if there is a group $\hat{G}$ having $G$ as a subgroup and 
a tuple $\bar{g}$ from $\hat{G}$, such that every equation/inequation from $p(\bar{x})$ is satisfied by $\bar{g}$ in $\hat{G}$. 
 Let $S^G_{W(\bar{x})}$ be the set of all atomic types over $W(\bar{x})$ which are finitely realized by tuples from $G$.   
We consider $S^G_{W(\bar{x})}$ together with the ordering by inclusion. 
Note that the atomic type $p(\bar{x})$ in the definition below is maximal. 

\begin{definition} \label{overW} 
Let $G\le \hat{G}$ and $\hat{G} = \langle G, \bar{g} \rangle$, and let 
$p(\bar{x})$ be the atomic type of $\bar{g}$ in $\hat{G}$ over $W(\bar{x})$: 
the set of all  equations/inequations defined for words from $W(\bar{x})$, which are satisfied by $\bar{g}$. 
We say that $\hat{G}$ is a limit group over $G$ with respect to $W(\bar{x})$ if 
every finite subset of $p(\bar{x})$ is realized (=satisfied) in $G$ by some tuple. 
\end{definition} 

\noindent 
It is worth noting that when $W(\bar{x}) = \mathbb{F}_{\bar{x}}$ and $G$ is finitely generated, 
this definition gives limit groups over $G$: fixing an appropriate tuple of markers, 
the group $\hat{G}$ appears as a limit satisfying Definition \ref{dlg}. 

When $W(\bar{x})$ is countable, the group $\hat{G}$ can be viewed as a limit of a sequence 
$(G, \bar{g}_n )$, $n\in \mathbb{N}$, where $\bar{g}_n$ is a realization in $G$ of first $n$ relations 
 of $p(\bar{x})$ corresponding to some fixed enumeration of $W(\bar{x})$. 
In this case we will write: $(\hat{G}, \bar{g} ) = \lim_{W(\bar{x})} (G, \bar{g}_n )$. 

\begin{remark} \label{compact} 
{\em Note that when a set $q(\bar{x})$ of relations of the form $w(\bar{x}) = 1$ or $w(\bar{x}) \not=1$ where 
$w(\bar{x}) \in W(\bar{x})$, is logically consistent
(with respect to $Th(G)$), 
then it is an atomic type from 
$S^G_{W(\bar{x})}$, and furthermore,  
there is a group $\hat{G}$ and a tuple $\bar{g} \in \hat{G}$ such that  
$\hat{G}$ is a limit group $\langle G, \bar{g} \rangle$ with $\bar{g}$  
realizing $q(\bar{x})$. 
Indeed, by the compactness theorem from model theory, there is a group $\hat{\mathsf{G}}$ 
\footnote{an elementary extension of $G$}  
and $\bar{g}$ in it that realizes $q(\bar{x})$.   
Then let $\hat{G}$ be the subgroup of $\hat{\mathsf{G}}$ generated by $G$ and $\bar{g}$. }
\end{remark} 

\noindent 
We emphasize here that $W(\bar{x})$ can be a proper subset of $G* \mathbb{F}_{\bar{x}}$. 
In particular, even if $\hat{G}$ is finitely generated, it is not necessary a $G$-limit in the sense of Definition \ref{dlg}. 
This is why we use the term {\em relative limits}. 
This issue will be treated in Section 4. 

The following easy proposition shows why free products are natural candidates for ``real" limits. 

\begin{proposition}\label{pl1} 
Let $G$ be a group and $W(\bar{x})$ be a set of non-trivial words from $G*\mathbb{F} _{m}$. 
Let $g_1 ,\ldots , g_m$ be a free basis of a free group of rank $m$. 
Then the existence of a limit group over $G$ that satisfies all the inequalities $\{ w(\bar{x}) \not= 1 \, | \, w(\bar{x}) \in W(\bar{x})\}$  
implies that $G*\langle g_1, \ldots , g_m \rangle$ is a limit group over $G$ with respect to $W(\bar{x})$. 
\end{proposition}

Mixed identities naturally arise in this context. 

\begin{definition} \label{law}
 Let $w(y_{1},\ldots , y_{\ell})$ be a non-trivial word over a group $G$, reduced in $G* \mathbb{F} _{\ell}$ 
and containing at least one variable. 
We will call $w(\bar{y})$ a mixed identity (or  law with constants) in $G$ if for every 
$\bar{g}=(g_{1},\ldots , g_{\ell})\in G^{\ell}$, 
 $w(\bar{g})=1$.
\end{definition}

A group $G$ is {\em mixed identity free} (shortly MIF) if it does not have mixed identities. 
By Proposition 5.3 from \cite{HO} this exactly means that $G$ and every $G*\mathbb{F}_n$ have the same universal theory. 
This implies the following proposition. 

\begin{proposition} \label{pmif} 
Assume that $G$ is MIF. 
Then $G*\mathbb{F}_n$ is a limit group over $G$. 
\end{proposition} 

Our main objects below are not MIF. 
In Section 3 we will show that they satisfy strikingly different properties. 
However, in Sections 2 and 4 we will show that our previous paper \cite{Iva-Za} gives a method for (relative) results resembling Proposition \ref{pmif}. 

\subsection{Thompson's group F} 

 To introduce the Thompson's group $F$ we will follow \cite{CFP}. 
{\em Thompson's group} $F$ is the group given by the following infinite group presentation:
\[ 
\langle \mathsf{x}_0, \mathsf{x_1, x}_{2}, \ldots \Big| \, \mathsf{x_{j}x_{i}=x_{i}x}_{j+1} , \, i<j  \rangle . 
\] 
In fact $F$ is finitely presented:
\[ 
F = \langle \mathsf{x}_0, \mathsf{x}_1\ \Big|\ [\mathsf{x}_0\mathsf{x}_1^{-1}, \mathsf{x}_0^{-i}\mathsf{x}_1\mathsf{x}_0^{i}] \, , \, i=1,2 \rangle .  
\] 

We will use the following geometric interpretation of $F$. 
Consider the set of all strictly increasing continuous piecewise-linear functions from the closed unit interval onto itself. 
Then the group $F$ is realized by the set of all such functions, which are differentiable except at finitely many dyadic rational numbers and such that all slopes (deriviatives) are integer powers of 2. 
The corresponding group operation is just the composition. 
For further reference it will be useful to give an explicit form of the generators $\mathsf{x}_{n}$, for $n\geq 0$, in terms of piecewise-linear functions:
\[ 
\mathsf{x}_{n}(t) = \left\{ \begin{array}{ll} t & \textrm{, $t\in [0,\frac{2^{n}-1}{2^{n}} ]$} \\
	\frac{1}{2}t + \frac{2^{n}-1}{2^{n+1}} & \textrm{, $t\in [\frac{2^{n}-1}{2^{n}}, \frac{2^{n+1}-1}
	{2^{n+1}} ]$} \\ t - \frac{1}{2^{n+2}} & \textrm{, $t\in [\frac{2^{n+1} -1}{2^{n+1}},
	\frac{2^{n+2}-1}{2^{n+2}}]$} \\	2t-1 & \textrm{, $t\in [\frac{2^{n+2}-1}{2^{n+2}},1] .$} \end{array}\right. 	 
\] 
For an arbitrary element $g$ in $F$ let $\mathsf{supp}(g)$ be the set $\{ t\in [0,1] : g(t) \neq t \}$
and $\overline{\mathsf{supp}}(g)$ the topological closure of $\mathsf{supp}(g)$. 

For any dyadic subinterval $[a,b]\subset [0,1]$, consider the set of elements in $F$, 
which are trivial on the complement, and denote it by $F _{[a,b]}$. 
We know that it forms a subgroup of $F$, which is isomorphic to the whole group. 
Let us denote its standard infinite set of generators by $\mathsf{x}_{[a,b],0}, \mathsf{x}_{[a,b],1}, \mathsf{x}_{[a,b],2}, \ldots$, 
where for $n\geq 0$ we have:
\[ 
\mathsf{x}_{[a,b],n}(t) = \left\{ \begin{array}{ll} t & \textrm{, $t\in [0,a+\frac{(2^{n}-1)(b-a)}{2^{n}}]$} \\
	\frac{1}{2}t + \frac{1}{2}(a+\frac{2^{n}-1}{2^{n}}) & \textrm{, $t\in [a+\frac{(2^{n}-1)(b-a)}{2^{n}},
			a+\frac{(2^{n+1}-1)(b-a)}{2^{n+1}}]$} \\
	 t - \frac{b-a}{2^{n+2}} & \textrm{, $t\in [a+\frac{(2^{n+1}-1)(b-a)}{2^{n+1}},
			a+\frac{(2^{n+2}-1)(b-a)}{2^{n+2}}]$} \\
	2t - b & \textrm{, $t\in [a+\frac{(2^{n+2}-1)(b-a)}{2^{n+2}},b]$} \\
	t & \textrm{, $t\in [b,1] .$} \end{array}\right. 
\] 
Moreover, if $\iota _{[a,b]}$ denotes the natural isomorphism between $F$ and $F_{[a,b]}$ sending $\mathsf{x}_{n}$ 
to $\mathsf{x}_{[a,b],n}$ for all $n\geq 0$, then for any $f\in F$ by $f_{[a,b]}$ we denote the element 
$\iota _{[a,b]} (f)\in F_{[a,b]}<F$. 

We will repeatedly use the following fact (see Lemmas 4.2 and 2.4 in 
\cite{CFP} and  \cite{KM}, respectively). 

\begin{fact} \label{CFP} 
If $0=s_0< s_1<s_{2}<\ldots <s_{n}=1$ and 
$0=t_{0}<t_{1}<t_{2}<\ldots <t_{n}=1$ are partitions of $[0,1]$ consisting of dyadic rational numbers, then there exists $f\in F$ such that $f(s_{i})=t_{i}$ for $i=0,\ldots , n$. 
	
Furthermore, if $s_{i-1}=t_{i-1}$
and $s_{i}=t_{i}$ for some $i$ with $1\leq i\leq n$, then $f$ can be taken to be the identity on the interval $[s_{i-1},s_{i}]$.
\end{fact}

Many statements and examples concerning Thompson's group, can be easily exposed using 
the \emph{rectangle diagrams} introduced by W. Thurston. 
The paper \cite{Zarphd} contains some pictures that can be helpful below.  

Since $F$ is a permutation group (on $[0,1]$), we sometimes apply standard terminology
and notation of the area of permutation groups (see \cite{H}). 
For any group $G$ acting on some set $X$ and any subset $A\subseteq X$, by $\mathsf{stab}_{G}(A)$
we denote the set of all elements from $G$, which stabilize $A$ pointwise. 

It is proved in \cite{Zar} and \cite{SZ} that $F$ is not MIF. 
As we mentioned in the end of Section 1.1, this has strong consequences 
concerning opposition of $F$-limits to $\mathbb{F}_2$-limits. 
On the other hand, as we have mentioned above, Brin's theorem shows that 
$\mathbb{F}_2$-limit groups are $F$-limits too. 
Thus, it is not obvious that free constructions involving $F$ 
cannot be realized as $F$-limits. 
In the next section we give examples when it is possible. 

\begin{remark} \label{prod} 
{\em It is interesting to note that direct products of copies of $\mathbb{F}_2$ are not limits of free groups, see \cite{GC}. 
On the other hand, the direct product $F \times F$ is fully residually $F$, i.e. a limit with respect to the class of subgroups of $F$. 
This follows from the fact that $F\times F$ is embeddable into $F$.}  
\end{remark}

\section{Limit groups over hereditarily separating groups}

The main result of this section is as follows. 
Consider the following subgroup of Thompson's group $F$
\[ 
H_{<1}:=\Big\{ f\in F\ \Big|\ 1\notin\overline{\mathsf{supp}}(f)\Big\} .
\] 
Consider the centralized HNN extension of $F$ relative to $H_{<1}$: 
\[ 
 F\ast _{H_{<1}} =\langle F, g\, | \, ghg^{-1} =h , \, h\in H_{<1}\rangle . 
\]
In Section 2.3 we will show that it is an $F$-limit group. 

Since we strongly use the paper \cite{Iva-Za}, we review some material from it in Sections 2.1. 
 In Section 2.2 we give some other applications of this approach.

\subsection{Inequalities in groups with hereditarily separating action}

In \cite{Iva-Za} the authors have introduced hereditarily separating actions 
on topological spaces, a generalization of separating actions studied by Abert in \cite{A}. 
These objects have become convenient for applications of some theory of inequalities over groups.   
We now review this material.

\paragraph{(A) Hereditarily separating actions by homeomorphisms} 

\bigskip 

\begin{definition}
	Let $G$ be a permutation group acting on an infinite set $X$. \begin{itemize} 
\item We say that $G$ separates $X$ if, for every finite subset $Y\subset X$, the pointwise stabilizer 
$\mathsf{stab}_{G}(Y)$ does not stabilize any point outside $Y$. 
\item Assume that $\mathcal{X}$ is a Hausdorff topological space, $G$ consists of homeomorphisms of $\mathcal{X}$ and the set of fixed points, $\mathsf{Fix}(G)$, is finite. 
We say that $G$ hereditarily separates $\mathcal{X}$ if, for every  open and infinite subset $Z\subseteq \mathcal{X}$ and for every finite subset $Y\subset Z$, the subgroup 
$\mathsf{stab}_{G}((\mathcal{X}\setminus Z)\cup Y)$ 
does not stabilize any point from
	$Z\setminus (Y\cup \mathsf{Fix}(G))$.
\end{itemize} 
\end{definition}

\begin{example} \label{alt} 
{\em It is pointed out in Section 2 of \cite{Iva-Za} that 
\begin{itemize} 
\item Thompson's group $F$ is hereditarily separating with respect to its standard action on $[0,1]$;  
\item  the action of any finitely generated weakly branch group on the boundary space of 
the corresponding infinite rooted tree is also hereditarily separating. 
\end{itemize} }
\end{example} 

\noindent 
The case of weakly branch groups will be discussed in Section 4.2, 
and the necessary definitions will be given there.

We now give examples that are not finitely generated. 
Consider $\mathbb{Q}$ as a topological space with respect to the standard topology induced from $\mathbb{R}$. 
Then $A(\mathbb{Q}) = \mathsf{Aut}(\mathbb{Q}, < )$, the group of order preserving permutations of rational numbers, acts on 
$\mathbb{Q}$ by homeomorphisms. 
Furthermore, let $S(x,y,z)$ be the ternary relation defining the circular order on $\mathbb{Q}$. 
Then $\mathsf{Aut} (\mathbb{Q}, S)$ acts on the circle $\mathbb{S}^1$ by homeomorphisms. 
The following proposition is straightforward. 

\begin{proposition} \label{Qos} 
The groups $A(\mathbb{Q})$ and $\mathsf{Aut} (\mathbb{Q},S)$ have hereditarily separating actions on the corresponding 
topological spaces on $\mathbb{Q}$. 
\end{proposition} 

\paragraph{(B) Inequalities} 

Let $G$ be a group. 
Any inequality over $G$ can be considered as follows. 
Let $w(\bar{y})$ be a word over $G$ on $\ell$ variables 
$y_{1},\ldots , y_{\ell}$. 
It is usually assumed that it is a nontrivial element of $\mathbb{F} _{\ell}\ast G$ where $\mathbb{F} _{\ell} = \mathbb{F}_{\bar{y}}$. 
In order to study the inequality $w(\bar{y}) \not= 1$, 
we take $w(\bar{y})$ in the reduced form. 
If $w(\bar{y})\notin\mathbb{F} _{\ell}$, we usually assume that it is of the form 
\begin{equation}  \label{1.1}  
w(\bar{y})=\mathsf{u}_{n}v_{n}\mathsf{u}_{n-1}v_{n-1}\ldots \mathsf{u}_{1}v_{1}, 
\end{equation}
where $n\in\mathbb{N}$, $\mathsf{u}_{i}\in \mathbb{F}_{\bar{y}}$ 
and $v_{i}\in G\setminus\{ 1\}$ for each $i\leq n$. 
It is clear that every word with constants is conjugate to a word in this form. 

\begin{definition} \label{prod-const} 
If in the form (\ref{1.1}),  $v_{n}\cdot \ldots \cdot v_{1}\neq 1$ then we say that the word $w(\bar{y})$ has non-trivial product of constants  
(in $\mathsf{supp}(v_{n}\cdot \ldots \cdot v_{1})\subseteq X$).
\end{definition} 

\noindent 
Then it is clear that the tuple of units $\bar{1}$ solves the inequality  $w(\bar{y}) \not=1$. 

Let now $G$ be a permutation group on $X$. 
We distinguish some specific types of words over $G$ with
respect to the action on $X$. 
Let $w(\bar{y})$ be in the form (\ref{1.1}). 
Define: 
\[ 
O_{w}:=\bigcap _{i=0} ^{n-1} v_{0}^{-1}v_{1}^{-1}\ldots v_{i}^{-1}( \mathsf{supp}(v_{i+1})) ,
\] 
where $v_{0}=\mathsf{id}$. 
If $w(\bar{y})\in\mathbb{F}_{\ell}$ then let 
$O_{w}:=X\setminus \mathsf{Fix}(G)$.

\begin{definition} \label{osc} 
Let $V\subseteq X$. 
We say that the word $w(\bar{y})\in\mathbb{F} _{\ell}\ast G$ 
is  explicitly oscillating in $V$ if $w(\bar{y})$ is non-trivial and $V\cap O_{w}\neq\emptyset$. 

When $V=X$ we just say that $w(\bar{y})$ is explicitly oscillating.
\end{definition} 

\noindent 
In the situation where $G$ acts on a Hausdorff topological space $\mathcal{X}$ by homeomorphisms the set $O_{w}$ is open. 

\noindent 
{\bf Notation.} 
We introduce the following notation. 
For an element $v\in G$ let $v^1 = v$ and $v^0 = \mathsf{id}$.   
For every $w(\bar{y})$ given in the form 
(\ref{1.1}) and for every set $A\subseteq X$ let 
\[ 
\mathcal{V}_{w}(A) = \{  v^{\varepsilon_j}_{j}\cdot \ldots \cdot v^{\varepsilon_1}_{1}(A)\, | \, (\varepsilon_1 ,\ldots , \varepsilon_j ) \in \{ 0,1 \}^j \, , \, 1 \le j \le n \} . 
\]

\paragraph{(C) Solving a system of oscillating inequalities}

We will use the notion of \emph{oscillating} words from \cite{Iva-Za}. 
The definition is too technical to be given in this paper. 
Intuitively, it describes words, that are explicitly oscillating after the transition from the region $O_w$ to some other place. 
In particular, explicitly oscillating words are oscillating. 
Furthermore, it is proved in Proposition 3.7 of \cite{Iva-Za} that words with non-trivial product of constants are also oscillating. 
This information is sufficient for us because the applications of \cite{Iva-Za} in our paper do not transcend these two cases. 
The following theorem is a simplified version of Theorem 3.8 from \cite{Iva-Za}. 

\begin{theorem} \label{uab}
Let $G$ act on a perfect metric space $(\mathcal{X},\rho )$ by
homeomorphisms. 
Let 
$\{ w_{1}(\bar{y}), w_{2}(\bar{y}), \ldots , w_{m}(\bar{y})\}$ 
be a set of reduced and non-constant words from  
$\mathbb{F} _{\ell}\ast G$ on $\ell$ variables  
$y_{1},\ldots ,	y_{\ell}$. 

If $G$ hereditarily separates $\mathcal{X}$ and every $w_{j}(\bar{y})$, $j\leq m$, is oscillating, then the set of inequalities 
$w_{1}(\bar{y})\neq \mathsf{id}, w_{2}(\bar{y})\neq \mathsf{id},\ldots , w_{m}(\bar{y})\neq \mathsf{id}$ has a solution in $G$. 

Moreover, in the case when all $w_i (\bar{y})$ are explicitly oscillating,  for any collection $\{ O_{j}\}$ such that $O_{j}$ is an open subset of the set $O_{w_{j}}$, $j\leq m$, there is a solution $(g_{1},\ldots , g_{\ell})$ of this set of inequalities such that 
$\mathsf{supp}(g_{i})\subseteq\bigcup _{j=1}^{m} (\bigcup \mathcal{V}_{w_{j}}(O_{j}))$ for $1\leq i\leq \ell$.
\end{theorem}

\subsection{Free subgroups of limit groups with respect to hereditary separating actions}

We now give an easy example of an application of Theorem \ref{uab}.

\begin{proposition}
Let $G$ be a finitely generated group having a hereditarily separating action on some perfect metric space $\mathcal{X}$ by homeomorphisms.
Suppose $h\in G$ such that $|h| = \infty$ and for every $\ell \in \mathbb{Z}\setminus \{ 0 \}$, $\mathsf{supp}(h^{\ell}) = \mathsf{supp}(h)$. 
Then there is a limit group $\hat{G}$ over $G$ of the form $\langle G,g \rangle$ such that $\{ h,g\}$ is a free basis of a free subgroup of $\hat{G}$. 
\end{proposition}
{\em Proof.} 
Assume $G = \langle g_{1},\ldots , g_{n}\rangle$. 
It follows from Definition \ref{osc} that each  non-trivial, non-constant word from 
$\langle h \rangle * \mathbb{Z}$ depending on $y$, is oscillating.
Let $W(y)$ be a finite set of non-trivial non-constant words from $\langle h \rangle * \mathbb{Z}$. 
By Theorem \ref{uab} there is $g'\in G$ realizing each inequality $w(y) \not= \mathsf{id}$ with $w(y) \in W(y)$.  
This defines a marked group $G_{W(y)} = \langle g_{1},\ldots , g_{n}, g' \rangle$. 
Applying this argument to all finite subsets of  non-trivial, non-constant words from $\langle h \rangle * \mathbb{Z}$ 
we get the statement of the proposition by compactness of $\mathcal{G}_{n+1}$. 
$\Box$ 

\bigskip 

This proposition explains how free groups appear in limits of $F$. 
Indeed, $F$ is torsion free and every non-trivial $h\in F$ satisfies the conditions above. 

We can develop this idea as follows. 
Let $W_{0}(\bar{y})$ be the set of words over $F$ with $\ell$ variables, 
which are reduced, non-constant in $\mathbb{F} _{\ell}\ast F$ and do not have constants 
outside of $\langle \mathsf{x}_0\rangle$; 
similarly, let $W_{1}(\bar{y})$ be the set of all such words with no constants outside of  $\langle \mathsf{x}_1\rangle$. 
Now for any $n\in\mathbb{N}$ let $(g_{n,1},\ldots ,g_{n,\ell})$ be a solution of a system of
inequalities $w_{1}(\bar{y})\neq \mathsf{id}, \dots , w_{n}(\bar{y}) \neq \mathsf{id}$, 
where $\{ w_{1}(\bar{y}),\ldots , w_{n}(\bar{y}),\ldots \}$ is an enumeration of $W_{0}(\bar{y})\cup W_{1}(\bar{y})$.
We consider a sequence $(G_{n})_{n<\omega}$, where
\[ 
G_{n}:=\Big(\langle \mathsf{x}_0, \mathsf{x}_1, g_{n,1},\ldots g_{n,\ell }\rangle , ( \mathsf{x}_0, \mathsf{x}_1, g_{n,1},\ldots ,
	g_{n,\ell} )\Big) . 
\] 
Similarly as above, we see that any limit group of the form $\lim G_{k_i}$ marked by a tuple
$(\mathsf{x}_0, \mathsf{x}_1,g_{1},\ldots , g_{t})$ contains both, 
$\langle \mathsf{x}_0\rangle\ast\langle g_{1},\ldots , g_{t}\rangle$
and $\langle \mathsf{x}_1\rangle\ast\langle g_{1},\ldots , g_{t}\rangle$.

\begin{remark} 
{\em It is worth noting here that Brin finds in \cite{Bri} a sequence of pairs $(g_i , h_i )$ generating $F$ such that 
$\lim (F, (g_i ,h_i )) =( \mathbb{F}_2, (x_1 ,x_2 ))$.  
This is a more complicated property than those given above. }
\end{remark}

There is a version of these results where finite generation is not assumed. 

\begin{proposition}
Let $G$ be a group having a hereditarily separating action on some perfect metric space 
$\mathcal{X}$ by homeomorphisms.
Suppose $h\in G$ such that $|h| = \infty$ and for every $\ell \in \mathbb{Z}\setminus \{ 0 \}$, $\mathsf{supp}(h^{\ell}) = \mathsf{supp}(h)$. 
Then there is a limit group $\hat{G}$ over $G$ with respect to the set of all inequalities $w(y) \not= \mathsf{id}$  
with non-trivial, non-constant words $w(y)$ from $\langle h \rangle * \mathbb{Z}$. 
In particular, $\hat{G} = \langle G,g \rangle$ such that $\{ h,g\}$ is a free basis of a free subgroup. 
\end{proposition}

\noindent 
{\em Proof.} 
The proof is the same as in the finitely generated case. 
At the final stage of it apply Remark \ref{compact} to get a limit group over $G$ with respect to the set 

$\{ w(y) \not= \mathsf{id} \, | \, w(y)$ is non-trivial and non-constant in $\langle h \rangle * \mathbb{Z}\}$. 
$\Box$ 

\bigskip 

\noindent 
In the special case, when $G = A(\mathbb{Q})$, we obtain that for any non-trivial $h\in G$ 
there is a limit group $\hat{G} = \langle G,g\rangle$ such that $\{ h,g \}$ is a free basis of  $\mathbb{F}_{2}$.

\subsection{HNN extensions of $F$ as $F$-limit groups}

In this section, we consider Thompson's group $F$ again. 
In the introduction to Section 2 we announced the following application of Theorem \ref{uab}.
Let $H_{<1}$ be defined as follows: 
\[ 
H_{<1}:=\Big\{ f\in F\ \Big|\ 1\notin\overline{\mathsf{supp}}(f)\Big\} .
\] 
Consider the centralized HNN extension of $F$ relative to $H_{<1}$: 
\[ 
 F\ast _{H_{<1}} =\langle F, g\, | \, ghg^{-1} =h , \, h\in H_{<1}\rangle . 
\] 
Note that $H_{<1}$ is not a finitely generated subgroup. 
We will later show (see Theorem \ref{mt}), it is a necessary condition for the property stated in the following theorem. 

\begin{theorem}\label{h1}
There is a convergent sequence of marked groups  \\ 
$G_{n} = (F, (g_{n}, \mathsf{x}_0, \mathsf{x}_1))$, $g_{n}\in F$, $n\in\omega$,  such that
\[ 
\lim_{n\to\infty}(G_{n})=( F\ast _{H_{<1}}, ( g,\mathsf{x}_0, \mathsf{x}_1)) . 
\] 
\end{theorem}

{\em Proof.} 
Let $W_{cf} (y)$ be the set of all non-trivial words from $F\ast _{H_{<1}} $, where $y$ corresponds to the stable letter. 
We want to show that 
\begin{quote} 
for every finite $W' \subset W_{cf} (y)$ and every finite $E\subset H_{<1} $ there is 
a non-trivial  element $g\in F$ such that $g$ commutes with each $h\in E$ and satisfies inequality $w(y) \not= \mathsf{id}$ 
for every $w(y) \in W'$. 
\end{quote} 
We may assume that every $w(y)\in W'$ is in the reduced form, i.e.  
$w(y)=y^{a_{k}}v_{k}\ldots y^{a_{1}}v_{1}$ where for all $i\leq k$ we have $a_{i}\in\mathbb{Z} \setminus\{ 0\}$ and 
$v_{i}\in F$ with $v_{i} \not\in H_{<1}$ for every $i>1$. 
We may also assume that when $v_1$ appears in $H_{<1}$, then $v_1 \in E$. 

Using the latter assumption, we obtain an equivalent task after all possible transformations of the following form: 
\[ 
w(y)=y^{a_{k}}v_{k}\ldots y^{a_2}  v_2 y^{a_{1}}v_{1}  \, \mbox{ with } v_1 \in H_{<1} \, \mbox{ and } \, v_2 \not= \mathsf{id} , \,  \mbox{ is replaced by }  
\] 
\[ 
\hspace{7cm} w'(y)=y^{a_1 + a_{k}}v_{k}\ldots y^{a_2} v_2v_1 . 
\] 
Note that when $w(y)$ has non-trivial $v_2$, no constants of $w'(y)$ belong to  $H_{<1}$. 
Thus, the support of each constant of such $w'(y)$ is cofinal in $[0,1)$. 
In particular, $w'(y)$ is explicitly oscillating, and $1 \in \overline{O_{w'}}$. 

If $w(y) =  y^{a_{1}}$ or $w(y) =  y^{a_{1}}v_{1}$, 
then it is explicitly oscillating by Definition \ref{osc}. 
Note that in the former case and when $v_1 \not\in H_{<1}$, again $1 \in  \overline{O_{w}}$. 
The latter case with $v_1 \in H_{<1}$ is different and will be treated separately. 
Let $W''$ be the part of $W'$ after removing all  words of the form  $w(y) =  y^{a_{1}}v_{1}$ with $v_1 \in H_{<1}$.  
As a result, we may assume that 
\begin{quote}
\hspace{2cm} each word of $W''$ is explicitly oscillating with 
$1 \in  \overline{O_{w}}$. 
\end{quote} 
Since all elements from $H_{<1}$ are trivial in some neighborhood of $1$ there is a number 
$n_0\in\mathbb{N}$ such that $\mathsf{supp}(h )\subseteq [0,1-\frac{1}{n_0}]$ for every  
$h \in E$. 
Enlarging $n_0$ if necessary, we satisfy 
\[ 
\hspace{-2cm} [1-\frac{1}{n_0},1)\subseteq\bigcap \{ O_{w} \, | \,  w(y) \mbox{ is a word from }  W'' \} \mbox{, and }  
\] 
\hspace{1,5cm} $\mathsf{supp}(h)\cap (\bigcup \mathcal{V}_{w}([1-\frac{1}{n_0},1) )= \emptyset$ for every $h\in E$ and every $w(y) \in W''$. 

Applying Theorem \ref{uab} to the set of words $W''$, we obtain a solution $g$ of the set of the corresponding inequalities 
$w(y)\neq\mathsf{id}$.
It also follows from Theorem \ref{uab} that we may choose a solution $g$ such that it is not trivial and 
$\mathsf{supp}(g)\cap \mathsf{supp} (h) = \emptyset$ for every $h\in E$. 
As a result $[g,h] = \mathsf{id}$ for every $h\in E$. 
Note that when $w(y) \in W'$ and $w(y) =  y^{a_{1}}v_{1}$ with $v_1 \in H_{<1}$ (i.e., in $E$) the support of $w(g)$ includes $\mathsf{supp}(g)$, i.e., $w(g) \not=\mathsf{id}$. 
This solves the task formulated in the beginning of the proof. 

To finish the proof of the theorem, take $\vec{H}_{<1 }=\{ h_{1}, h_{2}, \ldots \}$, an enumeration of $H_{<1}$. 
Let also $\vec{W}_{cf}=\{ w_1 (y) , w_{2}(y) ,\ldots\} \subseteq (F\ast\mathbb{Z}) \setminus F$ 
be an enumeration of the set of reduced words from $F*_{H_{<1}}$ with one variable $y$. 
Applying the procedure described above to each $\{ h_{1}, h_{2}, \ldots ,h_n \}$ and $\{ w_1 (y) , w_{2}(y) ,\ldots , w_n (y) \}$ 
we obtain the corresponding $g_n \in F$.  
In this way, we construct the whole sequence $(g_{n})_{n\in \omega}$. 
Taking a subsequence if necessary, we make the sequence of marked groups 
$G_n = (F, (g_n, \mathsf{x_0 , x_1} ))$, $n\in \mathbb{N}$, convergent. 
The construction ensures that $\lim G_n=( F\ast _{H_{<1}}, ( g,\mathsf{x}_0, \mathsf{x}_1))$.  
$\Box$ 

\section{Laws and limits. Limits of $F$} 

We will discuss restrictions on the form of $F$-limits. 
The following theorem is the main results of Section 3.

\begin{theorem}\label{mt}
	Let  $G_{n} = (F, (g_{n},\mathsf{x}_0, \mathsf{x}_1))$, $g_{n}\in F$, $n\in \mathbb{N}$, be a convergent sequence of groups, 
 and let $(G,(g,\mathsf{x}_0, \mathsf{x}_1))$ be its limit. 
Then	$G$ is not an HNN extension of $F$ of the form 
$F* _{\alpha}=\langle F, g\ |\ ghg^{-1} =	\alpha (h), \, h\in H\rangle$, 
where $H$ is a finitely generated subgroup of $F$ and $\alpha$ is an embedding of $H$ into $F$.
\end{theorem}

\noindent 
Note that Theorem \ref{h1} shows that Theorem \ref{mt} cannot be generalized to the case of infinitely generated
subgroups. 

The main tool used in this section is the fact that $F$ is not MIF; see \cite{Zar} and \cite{SZ}.  
Before the case of HNN extensions we consider the property of decomposition into free products. 
This case is much easier and is a good illustration how laws with constants work in our context. 

We also mention the following helpful remark. 

\begin{remark} 
{\bf  One-variable laws versus multi-variable laws.}
{\em According to Remark 5.1 of \cite{HO} if a group $G$ is non-trivial, then for every $n\in \mathbb{N}$ 
there is an embedding $G * \mathbb{F}_n \to G * \langle x \rangle$ which induces the identity map on $G$. 
Thus, 
\begin{quote} 
$\bullet$ if a group $G$ satisfies a law with constants in $n$ variables, then $G$ also 
satisfies an one-variable law with constants. 
\end{quote} } 
\end{remark}

\subsection{Free products}
 
We firstly formulate a general statement.

\begin{theorem}\label{pro}
Let $G$ be a group, which satisfies a law with constants, but which is not of finite exponent. 
Let $\hat{G}$ be the limit over $G$ of a convergent sequence of marked groups 
$G_{n} = ( G, (g_{n,1},\ldots , g_{n,\ell}))$, 
where $g_{n,1},\ldots , g_{n,\ell}\in G$, $n\in\mathbb{N}$. 

Then $\hat{G}\neq G\ast K$ for any non-trivial $K<\hat{G}$ which is not of finite exponent.
\end{theorem}

\noindent 
{\em Proof.} 
To obtain a contradiction suppose that $\hat{G}=G\ast K$, where 
$K$ is not of finite exponent, and $\hat{G} = \langle G, g_{1},\ldots , g_{\ell} \rangle$ such that  
$(g_{1},\ldots , g_{\ell})$ is the limit of  $(g_{n,1},\ldots , g_{n,\ell})$ from the formulation. 
Suppose that $G$ satisfies a law with constants. 
As we already know,  $G$ satisfies an one-variable law with constants;  
denote it by $w(y)$.  
Since $G$ is of infinite exponent,  so is $\hat{G}$, and the word $w(y)$ has constants from $G$. 

Let $\hat{g}=u(\bar{h},\bar{g})$ be an element of $K\setminus\{ 1\}$, where $\bar{h} \in G$. 
We assume that $|\hat{g}|$ is greater than the length of $w(y)$. 
Obviously, 
\[ 
w(u(\bar{h}, g_{n,1},\ldots , g_{n,\ell}))=_{G} 1 \mbox{  for all } n\in\mathbb{N} .
\]  
By the definition of $G$-limit groups, we see that
$w(u(\bar{h},\bar{g}))=1$ in $\hat{G}$. 

On the other hand, the reduced form of $w(y)\in G*\mathbb{Z}$ 
has non-trivial occurrences of $y$ and elements of $G$. 
Since $|\hat{g}|$ is greater than the length of $w(y)$, 
the equality $w(\hat{g})=_{\hat{G}}1$ contradicts  the fact that $\hat{G}$ is the free product of $G$ and $K$. 
$\Box$ 

\bigskip 

The group $A(\mathbb{Q})$ is a nice example for $G$ in this theorem. 
Indeed, it is torsion-free (see \cite{holland}) and has a law with constants (see \cite{BodSnTh}).  
We conclude that when $\hat{G}$ is a limit over $A(\mathbb{Q})$ of a convergent sequence of marked groups, 
then $A(\mathbb{Q})$  cannot be a direct summand of $\hat{G}$ (note that $\hat{G}$ must be torsion-free).  

We now apply Theorem \ref{pro} to limits of Thompson's group $F$.

\begin{corollary}\label{free}
Suppose that we are given a convergent sequence of marked groups \\ 
$\{ (G_{n}, (\mathsf{x}_0, \mathsf{x}_1,	g_{n,1},\ldots , g_{n,\ell})) \, | \, n\in\mathbb{N} \}$, 
where $G_{n}=F$, $g_{n,1},\ldots , g_{n,\ell}\in F$,	$n\in\mathbb{N}$. \\ 
Denote by $\hat{G}$ its limit. 

Then $\hat{G}\neq F*K$ for any non-trivial $K< \hat{G}$.
\end{corollary}

\noindent 
{\em Proof.} 
By \cite{SZ} and \cite{Zar} there is a word $w(y)$, which is a law with constants in $F$. 
Apply Theorem \ref{pro} for $G=F$. 
$\Box$ 

\bigskip 

Theorem \ref{pro} can be applied to all examples collected in Section 3 of \cite{BodSnTh}. 
Among $A(\mathbb{Q})$ let us mention the group $Sym (\mathbb{N})$, 
the group $\mathsf{Aut}(\mathbb{Q},S)$ of homeomorphisms of $\mathbb{S}^1$, 
the automorphism group of the generic semi-linear ordering with joins and the automorphism groups of Wa\.{z}ewski dendrites. 
Each of them has a law with constants and is not of finite exponent. 
In Section 4 we will discuss the case of branch groups.

\subsection{HNN extensions}

It turns out that when in Theorem \ref{mt} we additionally assume that the HNN extensions are centralized, the condition that $H$ is finitely generated can be replaced by the following assumption:
\begin{equation} \label{diamondsuit}
     \exists h',h''\in H\ \bigg( \{ 0,1\}\subseteq\overline{\mathsf{supp}}(h')\cup\overline{\mathsf{supp}}(h'')\bigg) . 
\end{equation} 
 Therefore, we prove Theorem \ref{mt} in two steps. 
First, we use (\ref{diamondsuit}) and prove the
centralized case.

\begin{theorem}\label{dmt}
	Let  $G_{n} = (F, (g_{n}, \mathsf{x}_0, \mathsf{x}_1))$, $g_{n}\in F$, $n\in \mathbb{N}$, be a convergent sequence of groups, 
and let $(G,(g,\mathsf{x}_0, \mathsf{x}_1))$ be its limit. 
Then $G$ is not an HNN extension of F of the form $F\ast _{H}=\langle F, g\ |\ [g,h]=\mathsf{id}, h\in H\rangle$, 
where $H$ is a subgroup of $F$ satisfying (\ref{diamondsuit}).
\end{theorem}

\begin{remark} 
{\em Note that condition (\ref{diamondsuit}) follows 
from the inclusion $[0,1]\subseteq  \bigcup \{ \overline{\mathsf{supp}} (h) \, | \, h\in H\}$ for finitely generated $H$. } 
\end{remark}

\paragraph{$\bullet$ Some further preliminaries concerning Thompson's group F} 
Consider the action of $F$ on $[0,1]$. 
Let $g\in F$. 
We will call each point from the set 
$P_{g} = (\overline{\mathsf{supp}}(g)\setminus \mathsf{supp}(g))\cap\mathbb{Z} [\frac{1}{2}]$ 
a \emph{dividing point} of $g$. 
This set is finite and defines a finite subdivision 
$[0,1] = [p_{0}, p_{1}]\cup [p_{1}, p_{2}]\cup \ldots \cup [p_{n-1}, p_{n}]$ where $p_0 = 0$ and $p_n = 1$. 
Taking restrictions of $g$ to these intervals, we obtain a presentation 
$g = g_{1}g_{2}\cdot \ldots \cdot g_{n}$ with $g_{i}\in F_{[p_{i-1}, p_{i}]}$ for each $i\le n$. 
Excluding the trivial $g_i$ we get a decomposition $g= g_{i_1}\cdot \ldots \cdot g_{i_k}$ with $g_{i_s} \not= \mathsf{id}$. 
It is called the \emph{defragmentation} of $g$. 

\begin{fact}\emph{(Corollary 15.36 in \cite{GS}, Proposition 3.2 in \cite{KM})}\label{GSf}
	The centralizer of the element $g\in F$ is the direct product of finitely many cyclic
	groups and finitely many groups isomorphic to $F$.
{\em Furthermore, if  $g$ has the defragmentation $g=g_{i_1} \cdot \ldots \cdot g_{i_k}$, 
then the cyclic components of the decomposition of the centralizer are generated by some roots of the corresponding  elements 
$g_{i_1}, \ldots , g_{i_k}$. 
The components of this decomposition, which are isomorphic to $F$, are just the groups of the form
	$F_{[a,b]}$, where $[a,b]$ is one of the subintervals $[p_{i-1},p_{i}]\subseteq [0,1]$ 
	which is stabilized pointwise by $g$.}
\end{fact}

\noindent 
Now it is easy to see that if $f,f',f''\in F$ and $f''$ is taken from the defragmentation of $f'$ 
then the equality $[f,f']=\mathsf{id}$ implies $[f,f'']=\mathsf{id}$.
This will be helpful below. 
\begin{fact}\label{GSc} (\cite{CFP})\, 
	The center of $F$ is trivial. 
In particular, the center of any $F_{[a,b]}$ is trivial. 
\end{fact}

Note that when we view  the elements of $F$ as functions, the relations 
$[\mathsf{x}_0\mathsf{x}_1^{-1}, \mathsf{x}_0^{-i}\mathsf{x}_1\mathsf{x}_0^{i}]$ for $i = 1,2$,  
just express that two functions commute if they have mutually disjoint supports, except of finitely many points. 
Furthermore, these relations imply their versions for $i>2$. 
Based on $\mathsf{x}_0^{-i}\mathsf{x}_1\mathsf{x}_0^{i} =\mathsf{ x}_{i+1}$, 
we conclude that all relations of the form $[\mathsf{x}_0\mathsf{x}_1^{-1},\mathsf{ x}_{k}]$, $k>1$, hold in Thompson's group $F$. 
 
 The following remark is another illustration of this point of view. 
\begin{quote}
$\bullet$ For any $f,f'\in F$ we have $\mathsf{supp}(f^{f'})=(f')^{-1}(\mathsf{supp}(f))$.  \, \, \,
\end{quote} 

The following statement is Theorem 1 from \cite{SZ}. 

\begin{proposition}\label{lwc2m}
	 Consider the standard action of Thompson's group $F$ on $[0,1]$. 
     Suppose that we are given
	two disjoint open dyadic subintervals $I_{i}=(p_{i}, q_{i})\subseteq [0,1]$, $i=1, 2$,
	and assume that $p_{1}<p_{2}$. 
Fix any non-trivial $h_{1}\in F_{\bar{I}_{1}}$ and $h_{2}\in	F_{\bar{I}_{2}}$ and denote:
	$$w^{-}(y)=[h_{1}^{y},h_{2}]=y^{-1}h_{1}^{-1}yh_{2}^{-1}y^{-1}h_{1}yh_{2}$$
	$$w^{+}(y)=[h_{1}^{(y^{-1})},h_{2}]=yh_{1}^{-1}y^{-1}h_{2}^{-1}yh_{1}y^{-1}h_{2}$$
	Then the word $w(y):=[w^{-}(y),w^{+}(y)]$ is a law with constants in $F$.
\end{proposition}

\paragraph{$\bullet$ Proofs of Theorems \ref{mt} and \ref{dmt}} 

\noindent 
\emph{Proof of Theorem \ref{dmt}.}\ 
To obtain a contradiction, assume that $G=F\ast _{H}$ where $g$ is the conjugator from the presentation.  
We start with the case where $H=F$. 
Then every $h\in F$ commutes with $g$. 
This implies that for every $h\in F$ there exists $n_{h}\in\mathbb{N}$ such that for all $n>n_{h}$ the equality 
$ [g_{n},h] =\mathsf{id}$ holds in $G_n$. 

Fix two elements $h_{1}, h_{2}\in F$ without dividing points so that they do not have a common root, and  
$\overline{\mathsf{supp}}(h_{1})=\overline{s\mathsf{upp}}(h_{2})=[0,1]$. 
By Fact \ref{GSf} there are roots $\hat{h}_{1}$ and $\hat{h}_{2}$ of $h_{1}$ and $h_{2}$ 
respectively such that  $C_{F}(h_{1})=\langle\hat{h}_{1}\rangle$ and $C_{F}(h_{2})=\langle\hat{h}_{2}\rangle$. 
Since for all $n>\max\{ n_{h_{1}}, n_{h_{2}}\}$, $g_{n}\in\langle\hat{h}_{1}\rangle\cap\langle\hat{h}_{2}\rangle$, 
we see that for almost all $n\in \mathbb{N}$ either $g_n$ does not have dividing points or  $g_{n}=\mathsf{id}$. 
In particular, in the former case some non-trivial $g_n$ has cyclic centralizer including $h_1$ and $h_2$; 
a contradiction to the choice of them.    
In the case where $g_n = \mathsf{id}$, we have $G=F\ncong F\ast _{F}$, 
a contradiction to the initial assumption of the proof.  

Consider the case $H\neq F$. 
Let $h_{1}, h_{2}, \ldots$ be an enumeration of $H$. 
For any $i\geq 1$ denote by $h_{i,1},\ldots h_{i,l_{i}}$, $l_{i}\in\mathbb{N}\setminus\{ 0\}$, 
the elements of the defragmentation of $h_{i}$. 
The following claim is the main part of the proof of the theorem. 
\bigskip 

\noindent 
 \textbf{Claim.} There is a non-trivial dyadic interval $[p,q]\subseteq [0,1]$ such that $F_{[p,q]}\cap\langle h_{ij}\
|\ i\geq 1, j\leq l_{i}\rangle$ is trivial or cyclic.\\

\noindent 
\emph{Proof of Claim.}\ Assume the contrary, i.e. 
\begin{quote} 
$\bullet$ for every dyadic $[p,q]\subseteq [0,1]$ the group 
$F_{[p,q]}\cap \langle h_{ij} \, | \,  i\geq 1, j\leq l_{i} \rangle$ is non-trivial and  
not isomorphic to $\mathbb{Z}$. 
\end{quote} 
Note that it implies that there is no non-trivial dyadic interval 
$[p,q]\subseteq [0,1]\setminus\bigcup _{j=1}^{\infty}\overline{\mathsf{supp}}(h_{j})$. 
Thus, we may assume that
$[0,1]\setminus\bigcup _{j=1}^{\infty}\overline{\mathsf{supp}}(h_{j})$ consists of isolated points. 

{\bf Case 1.}  Suppose that there is some number $s$ such that the set 
$[0,1]\setminus\bigcup _{i=1} ^{s}\overline{\mathsf{supp}}(h_{i})$ does not contain an interval. 
This obviously  implies that $\bigcup _{i=1} ^{s}\overline{\mathsf{supp}}(h_{i}) = [0,1]$. 

Since for any $i\leq s$ and any $j\leq l_{i}$, the group
\[ 
F_{\overline{\mathsf{supp}}(h_{ij})}\cap\langle h_{i',j'}\, |\, i'\geq 1, j'\leq l_{i'}\rangle 
\]
is neither trivial nor cyclic, for every pair $(i,j)$ as above there is an element
\[ 
h_{ij}'\in F_{\overline{\mathsf{supp}}(h_{ij})}\cap\langle h_{i',j'}\, |\, i'\geq 1, j'\leq l_{i'}\rangle , 
\] 
that does not have a common root with $h_{ij}$. 
(It is worth noting here that all elements which have common roots with $h_{ij}$ 
belong to the centralizer of $h_{ij}$ and the latter is cyclic.)
Enlarging $s$ if necessary we may assume that there is a covering of $[0,1]$ by a family of  
intervals of the form $\overline{\mathsf{supp}}(h_{ij})$ with $i\le s$ such that 
for every interval from this covering, say with $i \leq s$ and $j\leq l_{i}$, 
the corresponding  $h'_{ij}$ already belongs to  
$F_{\overline{\mathsf{supp}}(h_{ij})}\cap\langle h_{i'j'}\, |\, i'\leq s, j\leq l_{i'}\rangle$. 

Since for every $i$ the limit group $G$ satisfies $gh_{i}g^{-1} = h_{i}$, there is a natural number $n_0$ 
such that $F$ satisfies $g_{n}h_{i}g_{n}^{-1}= h_{i}$ for all $i\leq s$ and $n>n_0$.
This implies that for every $i\leq s$ and every $j\leq l_{i}$ all equations of the form 
$g_{n}h_{ij}g_{n}^{-1}= h_{ij}$ (and the corresponding $g_{n}h'_{i,j}g_{n}^{-1}= h'_{i,j}$) are satisfied in $F$ for all $n>n_0$. 

We claim that the conclusion of the paragraph above implies $g_{n}=\mathsf{id}$ for all $n>n_0$
(contradicting to the assumption that $g_{n}$ converges to $g\neq\mathsf{id}$). 
To see this, fix $n>n_0$ and consider the defragmentation of $g_{n}$. 
If it is non-empty take any non-trivial element of this defragmentation, say $g_{n,1}$.  
Since $g_{n}$ commutes with each $h_{ij}$ from the covering above, 
the element $g_{n,1}$ stabilizes the ends of the corresponding  intervals $\overline{\mathsf{supp}}(h_{ij})$. 
This implies that there is a pair of elements 
$h_{i_{0}, j_{0}}$ and $h_{i_{0},j_{0}}'$ as above in $F_{\overline{\mathsf{supp}}(g_{n,1})}$, 
which do not have a common root and such that
\[ 
[g_{n,1},h_{i_{0}, j_{0}}]=[g_{n,1},h_{i_{0},j_{0}}']=\mathsf{id} . 
\]
Assuming that $g_{n,1} \not= \mathsf{id}$, it follows from Fact \ref{GSf} that
\[ 
h_{i_{0},j_{0}}  \, , \,  h'_{i_{0},j_{0}}  \in \,   C_{F_{\overline{\mathsf{supp}}}(g_{n,1})} ( g_{n,1} ) , 
\]
and $h_{i_{0},j_{0}}$ and $h'_{i_{0}, j_{0}}$ have a common root,  a contradiction.
This argument can be applied to any $g_{n,t}$. 
We see that $g_n$ is trivial. 

{\bf Case 2.} 
Now consider the case where for every $s\in\mathbb{N}$ 
the union $\bigcup _{i=1} ^{s}\overline{\mathsf{supp}}(h_{i})$ does not cover $[0,1]$. 
For each $i\geq 1$, $j\leq l_{i}$ denote by $[p_{i,j},q_{i,j}]$ the support of $h_{ij}$. 
It follows from the assumptions of the theorem that some initial and some final subintervals of $(0,1)$ are covered by the
supports of two elements from $H$. 
Thus, there is some non-trivial dyadic interval $[u,v]$, $u\neq 0$, $v\neq 1$, such that for every $s\in\mathbb{N}$,
\[ 
[u,v]\not\subseteq \bigcup _{i=1}^{s}\Big(\bigcup _{j=1}^{l_{i}} [p_{i,j},q_{i,j}]\Big) . 
\]
Since $[0,1]\setminus\bigcup _{j=1}^{\infty}\overline{\mathsf{supp}}(h_{j})$ consists of isolated points, 
we may assume that there are some two elements $h_{r}$ and $h_{t}$ such that 
$u\in\overline{\mathsf{supp}}(h_{r})$ and $v\in\overline{\mathsf{supp}}(h_{t})$. 
Let $h_{r,j_{1}}$ and $h_{t,j_{2}}$, where $j_{1}\leq l_{r}$, $j_{2}\leq l_{t}$, be two elements of 
the defragmentations of $h_r$ and $h_{t}$ respectively, such 
that $u\in\overline{\mathsf{supp}}(h_{r,j_{1}})$ and $v\in\overline{\mathsf{supp}}(h_{t,j_{2}})$.
Note that, since 
$[u,v]\nsubseteq\overline{\mathsf{supp}}(h_{r})\cup\overline{\mathsf{supp}}(h_{t})$, the interval 
$\overline{\mathsf{supp}} (h_{r,j_{1}})$ does not intersect $\overline{\mathsf{supp}}(h_{t,j_{2}})$. 

Using the fact that the groups
\[ 
\langle h_{ij}\, |\, i\geq 1, j\leq l_{i}\rangle\cap F_{[p_{r,j_{1}},q_{r,j_{1}}]}
\]
and
\[ 
\langle h_{ij}\, |\, i\geq 1, j\leq l_{i}\rangle\cap F_{[p_{t,j_{2}},q_{t,j_{2}}]}
\]
are neither trivial nor cyclic, one can show as in Case 1, that there is some $n_0$ 
such that for all $n>n_0$
\[ 
g_{n}|_{[p_{r,j_1},q_{r,j_1}]\cup [p_{t,j_2},q_{t,j_2}]}=\mathsf{id}|_{[p_{r,j_1},q_{r,j_1}]\cup [p_{t,j_2},q_{t,j_2}]}. 
\]
Since $g_{n}$ is continuous and is increasing, we also have that 
$g_{n}([0,p_{r,j_{1}}])=[0,p_{r,j_{1}}]$, $g_{n}([q_{r,j_{1}},p_{t,j_{2}}])=[q_{r,j_{1}},p_{t,j_{2}}]$ and 
$g_{n}([q_{t,j_{2}},1])=[q_{t,j_{2}},1]$ for all $n>n_0$. 

 Let $f_{0}\in F_{[p_{r,j_{1}},q_{t,j_{2}}]}$ be any element satisfying $f_{0}(q_{r,j_{1}})=p_{t,j_{2}}$. 
Then 
\[ 
[u,v]\subseteq [p_{r,j_{1}},q_{t,j_{2}}]\subseteq ( [p_{r,j_{1}},q_{r,j_{1}}]\cup\overline{\mathsf{supp}}
	(f_{0})\cup [p_{t,j_{1}},q_{t,j_{2}}] ) . 
\] 
Since for every $s$, 
\[ 
[u,v]\not\subseteq \bigcup _{i=1}^{s} \Big(\bigcup _{j=1}^{l_{i}} [p_{i,j},q_{i,j}] \Big) , 
\]
the element $f_{0}$ does not belong to $H$. 

Consider the word $w_{1}(y) =[y,y^{f_{0}}]$, where $y$ denotes a variable. 
We claim that 
\begin{equation} \label{bull} 
\mbox{ for every } n>n_0, \, \, \, w_{1}(g_{n})=\mathsf{id}. 
\end{equation} 
Fix some $n>n_0$. 
For any $z\in [0,p_{r,j_{1}}]$ we have:
\[ 
w_{1}(g_{n})(z)=g_{n}^{-1}f_{0}^{-1}g_{n}^{-1}f_{0}g_{n}f_{0}^{-1}g_{n}f_{0}(z)=g_{n}^{-1}f_{0}^{-1}g_{n}^{-1}f_{0}g_{n}f_{0}^{-1}g_{n}(z). 
\]
However, $g_{n}(z)\in [0,p_{r,j_{1}}]$, and thus we obtain the further reductions of $w_{1}(g_n )$:
\[ 
w_{1}(g_{n})(z)=g_{n}^{-1}f_{0}^{-1}g_{n}^{-1}f_{0}g_{n}f_{0}^{-1}(g_{n}(z))=g_{n}^{-1}f_{0}^{-1}g_{n}^{-1}	(g_{n}^{2}(z))=z. 
\] 
Let $z\in [p_{r,j_{1}},q_{r,j_{1}}]$. 
Since $p_{t,j_{2}}=f_{0}(q_{r,j_{1}})$, we have $g_{n}(f_{0}(z))\leq f_{0}(q_{r,j_{1}})$. 
Hence $f_{0}^{-1}g_{n}f_{0}(z)\leq q_{r,j_{1}}$. 
On the other hand, since
\[ 
f_{0}([0,p_{r,j_{1}}])=[0,p_{r,j_{1}}]=g_{n}([0,p_{r,j_{1}}]) , 
\]
we also have $f_{0}^{-1}g_{n}f_{0}(z)\geq p_{r,j_{1}}$, i.e.,  
$f_{0}^{-1}g_{n}f_{0}(z)\in[p_{r,j_{1}},q_{r,j_{1}}]$. 
As a result, 
\[ 
w_{1}(g_{n})(z)=g_{n}^{-1}f_{0}^{-1}g_{n}^{-1}f_{0}g_{n}(f_{0}^{-1}g_{n}f_{0}(z))=g_{n}^{-1}f_{0}^{-1}g_{n}^{-1}f_{0}(f_{0}^{-1}g_{n}f_{0}(z))=g_{n}^{-1}(z). 
\]
Since $z\in [p_{r,j_{1}},q_{r,j_{1}}]$, then  $w_{1}(g_{n})(z)=g_{n}^{-1}(z)=z$. 

If $z\in [q_{r,j_{1}},q_{t,j_{2}}]$, then $f_{0}(z)\in [p_{t,j_{2}},q_{t,j_2}]$. 
Therefore, 
\[ 
w_{1}(g_{n})(z)=g_{n}^{-1}f_{0}^{-1}g_{n}^{-1}f_{0}g_{n}f_{0}^{-1}g_{n}f_{0}(z)=
g_{n}^{-1}f_{0}^{-1}g_{n}^{-1}f_{0}g_{n}f_{0}^{-1}f_{0}(z)= 
\] 
\[ 
=g_{n}^{-1}f_{0}^{-1}g_{n}^{-1}(f_{0}g_{n}(z))=g_{n}^{-1}f_{0}^{-1}(f_{0}g_{n}(z))=z. 
\]
When $z\in [q_{t,j_{2}},1]$, we apply the same argument as in the case  $z\in [0,p_{r,j_{1}}]$. 
As a result, we have the statement (\ref{bull}). 
This implies that the equality $w_{1}(g)=\mathsf{id}$ holds in $G$. 

As we already know, $f_0 \not\in H$. 
Thus, by Britton's Lemma $w_{1}(g)$ cannot be reduced in the
HNN extension $F\ast _{H}$. 
We obtain a contradiction with $w_{1}(g)=\mathsf{id}$. 
This finishes the proof of the claim. 

\bigskip 

Let us fix some non-trivial dyadic interval $[p,q]\subseteq [0,1]$ such that 
$F_{[p,q]}\cap\langle h_{i,j}\ |\ i\geq 1, j\leq l_{i}\rangle$ is trivial or cyclic. 
In the latter case, choose $h$ with $F_{[p,q]}\cap\langle h_{ij}\ |\ i\geq 1, j\leq l_{i}\rangle = \langle h \rangle$.  
If $h$ is non-trivial, then shortening the interval $[p,q]$
if necessary, we may assume that $\overline{\mathsf{supp}}(h)=[p,q]$. 

Now find two pairwise disjoint closed and dyadic subintervals $I_{i}=[a_{i},b_{i}]\subseteq [p,q]$, $i=1, 2$, $a_{1}<a_{2}$, 
and choose two non-trivial elements $f_{1}\in F_{I_{1}}$ and $f_{2}\in F_{I_{2}}$ without dividing points. 
Then it follows from Proposition \ref{lwc2m} that the word 
\[ 
w_{2}(y) :=[y^{-1}f_{1}^{-1}yf_{2}^{-1}y^{-1}f_{1}yf_{2}\, , \, yf_{1}^{-1}y^{-1}f_{2}^{-1}yf_{1}y^{-1}f_{2}] 
\] 
satisfies $w_{2}(f)=\mathsf{id}$ for all $f\in F$. 
This implies that for all $g_{n}$, $w_{2}(g_{n})=\mathsf{id}$,
and thus $w_{2}(g)$ is equal to $\mathsf{id}$ in $G$.  

 On the other hand $w_{2}(y)$ is non-trivial, non-constant and contains two types of constants $f_{1}^{\pm
1}$ and $f_{2}^{\pm 1}$. 
By Fact \ref{GSf} when $h\not=\mathsf{id}$ neither $f_{1}$ nor $f_{2}$ have a common root with $h$. 
Since $H<\langle h_{ij}\ |\ i\geq 1, j\leq l_{i}\rangle$, we have $f_{1},f_{2}\notin H$. 
Now it follows from Britton's Lemma on irreducible words in HNN extensions that the word 
\[ 
w_{2}(g) =[g^{-1}f_{1}^{-1}gf_{2}^{-1}g^{-1}f_{1}gf_{2},gf_{1}^{-1}g^{-1}f_{2}^{-1}gf_{1}g^{-1}f_{2}] 
\] 
is non-trivial in $G=F\ast_{H}$. 
This contradiction finishes the proof of the theorem. 
$\Box$ 

\bigskip 

Now we can proceed to Theorem \ref{mt}. 

\bigskip 

\noindent 
\emph{Proof of Theorem \ref{mt}.}  To obtain a contradiction, assume that $G=F\ast _{\alpha}$. 
Let $\{ h_{1},\ldots , h_{s}\}$ be a generating set of $H$ and $h_{ij}$, $i\le s$, $j\le l_i$, 
be elements of defragmentations of these generators. 

Firstly, consider the centralized case, i.e., $\alpha = \mathsf{id}$.
If $\bigcup _{i=1} ^{s}\overline{\mathsf{supp}}(h_{i}) = [0,1]$, then condition (\ref{diamondsuit}) holds, and we
apply Theorem \ref{dmt}. 
If the union $\bigcup _{i=1} ^{s}\overline{\mathsf{supp}}(h_{i})$ does not cover $[0,1]$, 
then the statement of the Claim from the proof of Theorem \ref{dmt} holds:  
there exists some non-trivial dyadic subinterval $[p,q]\subseteq [0,1]$ such that 
\[ 
F_{[p,q]}\cap\langle h_{ij}\ | \, h_{ij} \mbox{ are members of defragmentations of elements of } H \rangle =\{\mathsf{id}\}. 
\] 
Indeed, since each element of $H$ is expressed as a word of $h_{ij}$, $i\le s$, $j\le l_i$, the support of 
such an element is contained in $\bigcup \{ \mathsf{supp}(h_{ij}) \, | \, i\le s, \, j\le l_i\}$.  
Thus, any $[p,q] \subset [0,1]\setminus \bigcup _{i=1} ^{s}\overline{\mathsf{supp}}(h_{i})$ works. 
Now just apply the final part of  the proof of Theorem \ref{dmt}, i.e.  construct a suitable law with constants and 
find a contradiction with Britton's lemma.  

Consider the situation where $G=F\ast _{\alpha}$ and $\alpha$ is an arbitrary embedding. 
Let $(g_{n})_{n<\omega}$ be a sequence in $F$, such that for every $i\leq s$ 
the equation $g_{n}h_{i}g_{n}^{-1}=\alpha (h_{i})$ holds for almost all $n$. 
Fix an  element $f\in F$ such that for all $i\leq s$, $fh_{i}f^{-1}=\alpha (h_{i})$.  
We now apply the argument above (i.e., of Theorem \ref{dmt}) to the sequence $(f^{-1}g_{n})_{n<\omega}$. 
This sequence converges to some element $g'$ (in fact, $g'=f^{-1}g$ where $g$ is a conjugator from the presentation of $G$). 
Each $f^{-1}g_{n}$ commutes with every $h_{i}$, $i\leq s$, and therefore belongs to $C_G (H)$. 
To get a final contradiction, apply the statement of the theorem in the centralized case. 
$\Box$ 

\bigskip 

In the final remark of this section we comment on a possible connection with the results of the paper \cite{Os}. 
In order to  provide a general recipe for constructing finitely presented condensed groups, 
D. Osin considers in \cite{Os} 
for an $n$-generated group $G$, an injective continuous map $Sub(G) \to \mathcal{G}_{n+1}$ 
defined by the formula $H \to G*_{H}$. 
We do not know if there is a condensed group of the form $F*_{H}$. 
Theorem \ref{mt} shows that when it exists, it cannot be a limit over $F$.  

\section{Relative limits with respect to hereditarily separating groups}

Our original motivation for this paper was to describe the groups that arise as limits of finitely generated marked groups, and since the methods also apply for branch groups, they have been included in the paper (see Section 4.2). 
However, recently we realized that the restriction of finite generation is not necessary for the approach. 
Indeed, when we apply Definition \ref{overW}, infinitely generated (even uncountable) groups 
become available for our methods. 
For example, oligomorphic groups studied in \cite{BodSnTh} are natural objects for this. 
Furthermore, Definition \ref{overW}  introduces limits relative to some classes of words. 
Thus, one can look for partial versions of statements which are impossible 
in the absolute form. 
This issue is introduced in more detail in Section 4.1 below. 
In Section 4.3 we consider some oligomorphic groups that resemble Thompson's group $F$.

\subsection{Relative limit groups}

As we already know, groups with mixed identities usually do not appear as free summands in their limits. 
On the other hand, if we consider relative limits, the situation changes.  
The following proposition is a general observation of this kind. 
It uses terminology of Section 2.1. 

\begin{proposition}\label{pl} 
Let $G$ be a group having a hereditarily separating action on some perfect metric space $X$. 
Let $W(\bar{y})$ be the set of all non-trivial words from $G* \mathbb{F} _{m}$, 
which are oscillating. 
Let $e_1 ,\ldots , e_m$ be a free basis of a free group of rank $m$. 
Then $G*\langle e_1, \ldots , e_m \rangle$ is a limit group over $G$ with respect to $W(\bar{y})$. 

 In particular, also suppose that $G=\langle \mathsf{g}_1,\ldots , \mathsf{g}_{q}\rangle$. 
Then there exists a sequence of tuples 
$\{ (g_{n,1},\ldots , g_{n,m}) \, | \, n\in \mathbb{N} \}$  
from $G$ such that 
\[ 
G_{n}=(G, (\mathsf{g}_{1},\ldots , \mathsf{g}_{q}, g_{n,1},\ldots , g_{n,m})), n\in\mathbb{N}, 
\] 
is a convergent sequence of marked groups such that for every $w(\bar{y})\in W(\bar{y})$ 
the inequality $w(g_{n,1},\ldots , g_{n,m})\neq 1$ is satisfied in almost all $G_{n}$.
\end{proposition}

\noindent 
{\em Proof.} 
Let $W(\bar{y})=\{ w_{\iota}(\bar{y}) \, | \, \iota \in I\}$. 
By Theorem \ref{uab} any system of inequalities 
$\{ w_{\iota_1}(\bar{y})\neq 1, w_{\iota_2}(\bar{y})\neq 1,\ldots , w_{\iota_n}(\bar{y})\neq 1\}$
has a solution in $G$. 
The first statement follows.  

The second statement follows from the first one. 
Indeed, enumerate $W(\bar{y}) = \{ w_i (\bar{y}) \, | \, i \in \omega \}$ 
and construct a sequence of marked groups $(G_{n})_{n<\omega}$, where 
$G_{n}= G$ with markers $\mathsf{g}_{1},\ldots , \mathsf{g}_{q}$, $g_{n,1},\ldots , g_{n,m}$ such that 
$(g_{n,1}, \ldots , g_{n,m})$ is a solution in $G$ of the system 
\[ 
\{ w_{1}(\bar{y})\neq 1, w_{2}(\bar{y})\neq 1,\ldots , w_{n}(\bar{y})\neq 1\}.
\]
By compactness, this sequence can be chosen to be convergent. 
$\Box$ 

\bigskip 

By this proposition and Proposition \ref{Qos} 
we see that the groups $A(\mathbb{Q}) *  \mathbb{F}_m$ and $\mathsf{Aut}(\mathbb{Q},S) *  \mathbb{F}_m$ are  limit groups over  
$A(\mathbb{Q})$ and resp. over $\mathsf{Aut}(\mathbb{Q},S)$ with respect to the set of all oscillating words with $m$ variables. 

\begin{corollary} \label{th-gr}
Let $W(\bar{y})$ be a set of words over Thompson's group $F= \langle \mathsf{x_1 , x_2} \rangle$ 
with $m$ variables that are oscillating with respect to the natural action of $F$ on the unit interval. 
Then there exists a sequence of marked groups 
\[ 
G_{n}:=\Big(\langle  \mathsf{x_0, x}_1 ,g_{n,1},\ldots , g_{n,m} \rangle , (  \mathsf{x_0, x}_1, g_{n,1},\ldots ,
		g_{n,m})\Big) , \, n \in \mathbb{N} 
\]
for $g_{n,i}\in F$, such that some of its limits over $F$ with respect to $W(\bar{y})$ is of the form $\langle g_{1},\ldots , g_{m}\rangle\ast F$ where $\langle g_{1},\ldots g_{m}\rangle\cong \mathbb{F}_m$. 
\end{corollary}

Since the set of oscillating words is large, this statement is quite strong. 
However, it is not clear which other groups can be realized instead $\mathbb{F}_m$ in the formulation. 
In the following proposition, we get new possibilities, but the set of words becomes smaller. 

\begin{proposition}\label{iva}
Let $F=\langle \mathsf{x}_0, \mathsf{x}_1\rangle$ 
be Thompson's group and let $H\le F$ and $H=\langle h_{1},\ldots, h_{\ell}\rangle$.  
Then the group $H *F$ is a limit group over $F$ with respect to the set of words  which have constants from $F$ and have non-trivial product of constants. 
\end{proposition}

\noindent 
{\em Proof.} 
We firstly show that any system of inequalities over $F$, 
$w_1(\bar{y})\neq\mathsf{id}, \ldots ,w_m(\bar{y})\neq\mathsf{id}$, of non-constant words on $\ell$ 
variables with non-trivial products of constants, 
has a solution $g_1 ,\ldots , g_{\ell}\in F$ such that $\langle g_1,\ldots, g_{\ell}\rangle\cong H$. We will apply the material of Section 2.1.  

 Fix $w_{1}(\bar{y}),\ldots ,w_{m}(\bar{y})$. 
Without loss of generality, we assume that words $w_j(\bar{y})$, $1\leq j\leq m$, 
are written in the following reduced and non-degenerate form:
\[ 
\mathsf{u}_{j,k_{j}}v_{j,k_{j}}\ldots \mathsf{u}_{j,2}v_{j,2}\mathsf{u}_{j,1}v_{j,1}, 
\]
where $\mathsf{u}_{j,i}$ depend only on the letters $y_{1},\ldots ,y_{\ell}$ and 
$v_{j,i}\in F\setminus\{\mathsf{id}\}$, $1\leq i\leq k_{j}$, where 
$k_{j}\in\mathbb{N}\setminus\{ 0\}$, $1\leq j\leq m$. 

 Since $w_{1}(\bar{y}),\ldots ,w_{m}(\bar{y})$ have non-trivial product of constants, we choose points \\ 
 $p_{1},\ldots p_{m}\in [0,1]$ such that $w_{j}(\mathsf{id},\ldots , \mathsf{id})(p_{j})\neq p_{j}$ for all $j$ with $1\leq j\leq m$. 
Next we choose some dyadic, non-degenerate, and closed interval
\[ 
U\subseteq  [0,1]\setminus \Big( \bigcup_{j=1}^{m} \{ w' (\mathsf{id},\ldots ,\mathsf{id})(p_{j}) \, | \, w' (\bar{y} ) \mbox{ is a final segment of } w_j (\bar{y} )\} \Big) . 
\] 
Define a tuple $\bar{g}=(g_{1},\ldots , g_{\ell})\in F^{\ell}$ as follows: $g_{i}=h_{i \, \, U}$, $1\leq i\leq \ell$. 
Here we use the notation from Section 1. 
Thus $H_{U}=\langle g_{1},\ldots , g_{\ell}\rangle<F_{U}$ and obviously is isomorphic to the original $H$. 
Moreover, for each $j$ with $1\leq j\leq m$,
\[ 
w_{j}(g_{1},\ldots , g_{\ell })(p_{j})=\mathsf{u}_{j,k_{j}}v_{j,k_{j}}\ldots \mathsf{u}_{j,2}v_{j,2}\mathsf{u}_{j,1}v
_{j,1}(g_{1},\ldots , g_{\ell })(p_{j})= 
\] 
\[ 
=v_{j,k_{j}}\ldots v_{j,2}v_{j,1}(g_{1},\ldots , g_{\ell })(p_{j})\neq p_{j} , 
\] 
i.e.  $\bar{g}$ is a solution of the set of inequalities $w_{1}(\bar{y})\not=\mathsf{id},\ldots ,w_{m}(\bar{y})\neq \mathsf{id}$, 
and generates a subgroup of $F$ isomorphic to $H$. 

 Finally, enumerate all non-constant words on $\ell$ variables with non-trivial product of constants: 
$w_1 (\bar{y}),\ldots , w_n (\bar{y}),\ldots$ . 
We construct a sequence of groups 
\[ 
G_n =(F, (g_{n,1},\ldots , g_{n,\ell }, \mathsf{x}_0 , \mathsf{x}_1 )), \, n\in \omega, 
\] 
where $g_{n,1},\ldots , g_{n,\ell }$ are defined as above with respect to 
$w_1 (\bar{y}),\ldots , w_n (\bar{y})$. \\ 
Since $\langle g_{n,1},\ldots , g_{n,\ell }\rangle\cong H$ for all $n$, 
we see that $G_H = \lim (G_n )_{n\in\omega}$ satisfies the statement of the
proposition. 
$\Box$ 

\bigskip   

There are some additional comments. 
Note that if $e_1, \ldots , e_{\ell}$ is a free basis of $\mathbb{F}_{\ell}$, then both Proposition \ref{pl} and Proposition \ref{iva} induce homomorphisms 
$\phi$ from $\mathbb{F}_{\ell} * F$ to some $F$-limit groups $\hat{F}>F$ 
such that each $\phi$ fixes $F$. 
In the latter case we also have a subgroup $H= \langle h_1 ,\ldots , h_{\ell} \rangle < F$ and the condition $\phi (\mathbb{F}_{\ell}) \cong H$. 
In each case a homomorphism $\psi: \phi (\mathbb{F}_{\ell}) \to H$ can be defined that 
takes every $\phi (e_i)$ to $h_i$. 
In the case of Proposition \ref{pl} this follows from the fact that $\phi$ is injective on $\mathbb{F}_{\ell}$.
In Proposition \ref{pl}, $\mathsf{Ker} \psi$ is large, but the kernel of $\phi(\mathbb{F}_{\ell}) * F \to \hat{F}$ 
is sufficiently small: all oscillating words over $F$ with $\ell$ variables are not in the kernel. 
On the other hand, Proposition \ref{iva} makes the kernel of $\psi$ trivial but the kernel from 
$\phi (\mathbb{F}_{\ell}) * F$ to $\hat{F}$ may increase.  
It only states that 
the words, which have non-trivial product of constants, are not in the kernel. 

This leads to the following research proposal. 
{\em Given 
$H =\langle h_1 ,\ldots , h_{\ell}\rangle <F$, investigate homomorphisms $\phi$ 
from $\mathbb{F}_{\ell} * F$ to limits $\hat{F}$ 
over $F$ such that the the map $\phi(e_i ) \to h_i$, $1 \le i \le \ell$, induces a homomorphism 
$\psi: \phi (\mathbb{F}_{\ell} )  \to H$.}  
Doing this task, pay special attention to the cases when $\mathsf{Ker} \psi$ and $\mathsf{Ker} \phi$ are minimal as possible.

\subsection{The case of branch groups} 
We shortly give some preliminaries on \emph{weakly branch groups}  
(following \cite{Gr, Gri} and \cite{A}). 
Let a group $G$ act isometrically on some rooted tree $T$ of finite valency. 
The vertices in $T$, which are at the same distance from the root are said to be at the same \emph{level}. 
We say that the action of the group $G$ on $T$ is \emph{spherically transitive} if $G$ acts transitively 
on each level of $T$. 
For any vertex $t\in T$ we define its 
{\em rigid stabilizer} to be the set of all elements from $G$, which stabilize  
$T\setminus T_t$ pointwise, where $T_t$ is the natural subtree hanging from $t$. 
A group $G$ is called  
{\em weakly branch} if it acts spherically transitively on some rooted tree $T$ so that the rigid stabilizer of every vertex is non-trivial.
 The \emph{boundary} of a tree $T$, denoted by $\partial T$, consists of the infinite branches
starting at the root. 
There is a topology on $\partial T$ where the base of open sets is determined by subtrees $T_t$. 
A weakly branch group obviously acts on the boundary $\partial T$ by homeomorphisms (in fact by isometries with respect to the natural metric). 
As we already know from Section 2.1 this action is hereditarily separated. 
Our first observation is as follows. 

\begin{proposition} 
Any weakly branch groups satisfies a mixed identity. 
\end{proposition}

\noindent 
{\em Proof.} 
Let $(G, \partial T)$ be a weakly branch group with the corresponding action on a tree. 
Let $t_1, t_2$ and $t_3$ be three distinct vertices of the same level. 
For every $i\le 3$ choose some non-trivial $h_i \in G$ from the rigid stabilizer of $t_i$.   
We claim that $G$ satisfies the following identity: 
\[ 
[[ y h_1 y^{-1}  , \, h_2 ], [yh_1 y^{-1} , \, h_3 ]] = \mathsf{id} . 
\] 
Indeed, if $\mathsf{supp}(y h_1 y^{-1}) \cap \mathsf{supp}(h_2) = \emptyset$, then  
$ [yh_1 y^{-1} , \, h_2 ] = \mathsf{id}$. 
If $\mathsf{supp}(y h_1 y^{-1}) \cap \mathsf{supp}(h_2) \not= \emptyset$, then 
$\mathsf{supp}(y h_1 y^{-1}) \cap \mathsf{supp}(h_3) = \emptyset$ and 
$ [yh_1 y^{-1} , \, h_3 ] = \mathsf{id}$. 
$\Box$ 

\bigskip 

Applying Theorem \ref{pro} we obtain the following corollary. 

\begin{corollary}\label{free}
Let $G$ be a weakly branch group. 	 
Suppose, we are given a convergent sequence of marked groups 
$((G_{n}, (g_{n,1},\ldots , g_{n,\ell}))\, | \,n\in \mathbb{N})$, 
where $G_{n}=G$, $g_{n,1},\ldots , g_{n,\ell}\in G$,	$n\in\mathbb{N}$. 
Denote by $\hat{G}$ its limit. 

Then $\hat{G}\neq G*K$ for any $K< \hat{G}$ of unbounded exponent.
\end{corollary}

The class of weakly branch groups contains the first Grigorchuk group (see \cite{Gr, Gri}), 
which is given by the following presentation:
\[ 
\langle \mathsf{a,b,c,d} \, | \, 1=\mathsf{a^{2}=b^{2}=c^{2}=bcd=\sigma ^{k}((ad)}^{4})=\sigma ^{k}((\mathsf{adacac})
	^{4}), k=0, 1,\ldots\rangle , 
\]
where the substitution $\sigma$ is defined by
\[ 
\sigma := \left\{ \begin{array}{ll} \mathsf{a}\to \mathsf{aca}\\
	\mathsf{b}\to \mathsf{d}\\
	\mathsf{c}\to \mathsf{b}\\
	\mathsf{d}\to \mathsf{c}.\end{array}\right. 
\]
We denote it by $\mathsf{G}$. 
Since $\mathsf{G}$ is weakly branch, Corollary \ref{free} can be applied to it. 
(Note that B. Jacobson has already observed in \cite{B-Jac} that $\mathsf{G}$ has mixed identities.) 
It is also known that $\mathsf{G}$ is a torsion group without infinite finitely generated subgroups of finite exponent. 
As a result we have the following statement. 
\begin{quote} 
$\bullet$ Suppose, we are given a convergent sequence of marked groups 
\[ 
(\mathsf{G}, (\mathsf{a}, \mathsf{b}, \mathsf{c}, \mathsf{d}, g_{n,1},\ldots , g_{n,\ell})\, | \, n\in \mathbb{N}), 
\] 
where $g_{n,1},\ldots , g_{n,\ell}\in \mathsf{G}$,	$n\in\mathbb{N}$. 
Denote by $\hat{G}$ its limit. 

Then $\hat{G}\neq \mathsf{G}*K$ for any infinite $K< \hat{G}$.
\end{quote} 
Let us also note that the statement of Corollary \ref{th-gr} holds in the case of this group too. 
\begin{quote} 
$\bullet$ Let $W(\bar{y})$ be a set of words over Grigorchuk's group  $\mathsf{G}$ with $m$ variables,  which are oscillating with respect to the natural action of $\mathsf{G}$ on an infinite rooted tree and on the boundary of it. 
Then there exists a sequence of marked groups 
\[ 
G_{n}:=\Big(\langle g_{n,1},\ldots , g_{n,m},  \mathsf{a}, \mathsf{b}, \mathsf{c}, \mathsf{d} \rangle , ( g_{n,1},\ldots , g_{n,m},  \mathsf{a}, \mathsf{b}, \mathsf{c}, \mathsf{d} )\Big)  
, \, g_{n,i}\in \mathsf{G}, \, n\in \mathbb{N} , 
\] 
such that $\langle g_{1},\ldots ,g_{m}\rangle\ast \mathsf{G}$ is its limit over $\mathsf{G}$ with respect to $W(\bar{y})$.
\end{quote} 
\bigskip 

Let $G$ be a weakly branch group acting on a rooted tree $T$. 
For $n\in\mathbb{N}$ let $T_n$  be  the subtree of $T$ consisting of the first $n$ levels. 
The group $G$ is called  \emph{self-reproducing} if for any $n$ and any $v$ of the $n$-th level, 
the pointwise stabilizer $\mathsf{stab}_G (T_n)$ induces on the subtree $T_{v}$ the action isomorphic to $(G,T)$; 
see page 142 of \cite{Gri} where this definition is given in a slightly weaker form. 
In the following proposition we show that the argument of Proposition \ref{iva} 
also works in the case of self-reproducing branch groups.

\begin{proposition}
 Let $G=\langle g_{1},\ldots , g_{n}\rangle$ be a self-reproducing weakly branch group with
respect to some action on a rooted tree $T$. 
Let $H\le G$ and $H=\langle h_{1},\ldots, h_{\ell}\rangle$. 
Then the group $G *H$ is a limit group over $G$ with respect to the set of words  which have constants from 
$G$ and have non-trivial product of constants. 
\end{proposition}

\noindent 
{\em Proof.} 
Similarly as in the proof of Proposition \ref{iva}, we show that any set of inequalities 
$w_1(\bar{y})\neq\mathsf{id}, \ldots ,w_m(\bar{y})\neq\mathsf{id}$, of non-constant words on $\ell$ 
variables with non-trivial products of constants, 
has a solution $g_1 ,\ldots , g_{\ell}\in G$ such that $\langle g_1,\ldots, g_{\ell}\rangle\cong H$.  

 Since $w_{1}(\bar{y}),\ldots w_{m}(\bar{y})$ have non-trivial product of constants, there are vertices \\ 
$t_{1},\ldots ,t_{m}\in T$ such that  for any $j\leq m$, the vertex $t_{j}$ is not stabilized by 
$w_{j}(\mathsf{id},\ldots , \mathsf{id})(x)$. 
Let $k$ the minimal level of the tree $T$ such that $t_{1},\ldots , t_{m}$ are contained within $T_k$. 

Let $k'>k$ and $t\in T$ belong to the $k'$-th level of $T$. 
It is easy to see that for any tuple of elements $\bar{f}=(f_{1},\ldots , f_{\ell}) \in G^{\ell}$ which fix $T_{k} \cup \{ t \}$,  
and any $j\leq m$ we have $w_{j}(\bar{f})(t_{j})\neq t_{j}$. 
Since $G$ is self-reproducing, the stabilizer $G_{t}$ is isomorphic to $G$,
and thus we may choose $\bar{f}$ so that $\langle f_{1},\ldots , f_{\ell}\rangle\cong H$ under the natural map. 

 We finish the proof exactly as in the case of Proposition \ref{iva}, 
i. e. enumerate all non-constant words on $\ell$ variables with non-trivial product of constants: 
$w_1 (\bar{y}),\ldots , w_n (\bar{y}),\ldots$. 
Construct a sequence of groups 
$G_n =(G, (f_{n,1},\ldots , f_{n,\ell } ))$, $n\in \mathbb{N}$, 
where $f_{n,1},\ldots , f_{n,\ell }$ is defined as above
with respect to $w_1 (\bar{y}),\ldots , w_n (\bar{y})$. 
Since $\langle f_{n,1},\ldots , f_{n,\ell }\rangle\cong H$ for all $n$, 
we see that $G_H = \lim (G_n )_{n\in\mathbb{N}}$ satisfies the statement of the
proposition. 
$\Box$ 

\bigskip 

 Since the first Grigorchuk group is an example of a self-reproducing weakly branch group (see \cite{Gri} for details), in this proposition $\mathsf{G}$ can be taken as $G$. 

\subsection{The case of order preserving permutations of $\mathbb{Q}$}

Now consider $A(\mathbb{Q})$. 
The following statement is a corollary of the proof of Proposition \ref{iva}. 

\begin{corollary}\label{iva-cor}
	 Let $G$ be $A(\mathbb{Q})$ or $\mathsf{Aut}(\mathbb{Q},S)$ and let $H\le A(\mathbb{Q})$ with 
$H=\langle h_{1},\ldots, h_{\ell}\rangle$. 
Then the group $G *H$ is a limit group over $G$ with respect to the set of words  which have constants from 
$A(\mathbb{Q})$ and have non-trivial product of constants. 
\end{corollary}

In the case of $A(\mathbb{Q})$ one can find a version of Theorem \ref{h1}. 
Indeed, let 
\[ 
\mathsf{O}_{cf} = \{ h\in A(\mathbb{Q}) \, | \, \mbox{ there is } q \in \mathbb{Q} \mbox{ such that } \forall x (x \ge q \to h(x) = x )\} 
\] 
be the stabilizer of the upper end of $\mathbb{Q}$. 

We say that a word 
\[ 
 w(y)=y^{a_{k}}v_{k}\ldots v_2 y^{a_{1}}v_{1}  
\] 
has {\em cofinal constants}, if every $v_i$ has an {\em orbital} that is cofinal in $\mathbb{Q}$. 
After reconstruction of the corresponding definitions from \cite{holland}, 
this condition just says that there is an element $q\in \mathbb{Q}$ such that for any $q'\in \mathbb{Q}$ some powers  
$v^{\ell}_i (q)$, where $\ell \in \mathbb{Z}$, are greater than $q'$. 

\begin{theorem} 
The group $(A(\mathbb{Q}) *_{\mathsf{O}_{cf}} \, , \, ( A(\mathbb{Q}),g))$ is a limit group over $A(\mathbb{Q})$ with respect to all words with cofinal constants. 
\end{theorem} 

\noindent 
{\em Proof.} 
The proof is similar to the proof of Theorem \ref{h1}. 
Let $W_{cf} (y)$ be the set of all non-trivial words from $A(\mathbb{Q})*_{\mathsf{O}_{cf}}$ with cofinal constants, 
where $y$ corresponds to the stable letter. 
We want to show that for every finite $W' \subset W_{cf} (y)$ and every finite $E\subset \mathsf{O}_{cf}$ there is 
a non-trivial  element $g\in A(\mathbb{Q})$ such that $g$ commutes with each $h\in E$ and satisfies the inequality 
$w(y) \not= \mathsf{id}$ for every $w(y) \in W'$. 
This task can be solved by a simplified version of the proof of Theorem \ref{h1}. 
$\Box$ 

\bigskip




\noindent 
{\bf Acknowledgements.}
The authors thank the referee for helpful comments. 


\section*{Declarations}

\begin{itemize}
\item Funding.  Not applicable.
\item Conflict of interest/Competing interests (check journal-specific guidelines for which heading to use). Not applicable.
\item Ethics approval and consent to participate. Not applicable. 
\item Consent for publication. Not applicable.
\item Data availability. No datasets were generated or analysed for this work.
\item Materials availability. Not applicable.
\item Code availability. Not applicable.
\item Author contribution. Not applicable.
\end{itemize}

\vspace*{10mm}

\begin{flushleft}
\begin{footnotesize}
Aleksander Iwanow

Institute of Computer Science, University of Opole, 

ul. Oleska 48, 45 - 052 Opole, Poland 

aleksander.iwanow@uni.opole.pl 

\bigskip 

Roland Zarzycki

Institute of Sociology, Wroc{\l}aw University, 

ul. Koszarowa 3,  50 - 350  Wroc{\l}aw, Poland 

bestializator@gmail.com

\end{footnotesize}
\end{flushleft}

\end{document}